\documentclass[12pt]{article}
\usepackage{amsmath}
\usepackage{amssymb}
\usepackage{braket} 
\usepackage{enumerate} 
\usepackage{amsthm}

\newtheorem{thm}{Theorem}
\newtheorem{prop}[thm]{Proposition}
\newtheorem{lem}[thm]{Lemma}

\numberwithin{equation}{section}
\theoremstyle{definition}
\newtheorem{dfn}[thm]{Definition}
\newtheorem*{ack}{Acknowledgement}

\newcommand{\Z}{\mathbb Z}
\newcommand{\IIi}{${\textrm{II}}_1$ }
\newcommand{\tr}{\mathrm{tr}}
\newcommand{\ce}[2]{\mathbb E_{#1}\left(#2\right)}
\newcommand{\nm}[1]{\left\|{#1}\right\|}
\newcommand{\VN}[1]{\mathcal L(#1)}
\newcommand{\vnotimes}{\overline{\otimes}}
\newcommand{\Norm}[1]{\mathcal N(#1)}
\newcommand{\Normo}[2]{\mathcal N^{#2}(#1)}
\newcommand{\NormN}[2]{\mathcal N_{#2}(#1)}
\newcommand{\Unit}[1]{\mathcal U(#1)}
\newcommand{\ip}[2]{\left<#1,#2\right>}
\newcommand{\HS}{\mathcal H}

\title{Semi-regular masas of transfinite length}

\author{Stuart White\thanks{\texttt{white@math.tamu.edu}}\and Alan Wiggins\thanks{\texttt{awiggins@math.tamu.edu}}}

\date{November 20, 2006}
\begin{document}
\maketitle
\begin{abstract}
In 1965 Tauer produced a countably infinite family of semi-regular masas in the hyperfinite \IIi factor, no pair of which are conjugate by an automorphism.  This was achieved by iterating the process of passing to the algebra generated by the normalisers and, for each $n\in\mathbb N$, finding masas for which this procedure terminates at the $n$-th stage.  Such masas are said to have length $n$.  In this paper we consider a transfinite version of this idea, giving rise to a notion of ordinal valued length.  We show that all countable ordinals arise as lengths of semi-regular masas in the hyperfinite \IIi factor.  Furthermore, building on work of Jones and Popa, we obtain all possible combinations of regular inclusions of irreducible subfactors in the normalising tower.
\end{abstract}

\section{Introduction}\label{Intro}
The study of maximal abelian self-adjoint subalgebras (masas) of \IIi factors dates back to Dixmier \cite{Dixmier.Masa}.  He considered the unitary normalisers of a masa $A$ in a \IIi factor $M$, that is the group $\Norm{A}$ (or $\NormN{A}{M}$ when the underlying factor $M$ is unclear) of all unitaries $u\in M$ with $uAu^*=A$.  The masa $A$ is said to be Cartan when $\Norm{A}''=M$ and singular when $\Norm{A}\subset A$.  In between we have the semi-regular masas; those where $\Norm{A}''$ is a proper subfactor of $M$.  Tauer constructed a countably infinite family of semi-regular masas in the hyperfinite \IIi factor such that no pair is conjugate by an automorphism \cite{Tauer.Masa}.  To show that the members of her family are non-conjugate, she introduced the length of a masa by iterating the normaliser construction and considering the chain below.
$$
A\subset\Norm{A}''\subset\Norm{\Norm{A}''}''\subset\dots
$$
The invariant Tauer considered was the length of this chain, i.e. how many iterations are required to reach the underlying factor $M$, if indeed this ever occurs.  In this paper we shall generalise this idea, considering masas of ordinal length.
\begin{dfn}
Let $B$ be a subalgebra of the \IIi factor $M$.  For ordinals $\alpha$ define $\Normo{B}{\alpha}$ inductively by
$$
\Normo{B}{\alpha}=\begin{cases}B&\alpha=0\\ \Norm{\Normo{B}{\beta}}''&\alpha=\beta+1\\ (\bigcup_{\beta<\alpha}\Normo{B}{\beta})''&\alpha\textrm{ is a limit ordinal}\end{cases}.
$$
We will say that $B$ has length $\alpha$ when $\alpha$ is the minimal ordinal for which $\Normo{B}{\alpha}=\Normo{B}{\alpha+1}$.
\end{dfn}
Our definition of length extends that of Tauer, who only considered the length of semi-regular masas $A$ in $R$ with $\Normo{A}{n}=R$ for some $n\in\mathbb N$.  In particular this definition ensures that if $M$ is faithfully represented on a Hilbert space $\HS$, then the length of every von Neumann subalgebra of $M$ has cardinality at most $\dim(\HS)$.

In \cite[section 6]{Tauer.Masa}, Tauer gives two non-conjugate length $2$ semi-regular masas in the hyperfinite \IIi factor $R$.  To distinguish them she introduced a property for subfactors $N\subset M$, which we shall call the \emph{Tauer property}: given any $x,y\in M$ with $\ce{N}{x}=\ce{N}{y}=0$, then $xy^*\in N$.  Her length $2$ semi-regular masas $A$ and $B$ in $R$ were constructed so that $\Normo{A}{1}\subset R$ has the Tauer property, while $\Normo{B}{1}\subset R$ does not.  Today, this property can be seen as an equivalent formulation of having Jones index $2$ and we will establish this in section \ref{Infinite}.

In \cite[section 2.4]{Saw.Thesis}, the first author constructed semi-regular masas in $R$ with length $\omega$, the first infinite ordinal, using matrix methods based on \cite{Tauer.Masa}.  Furthermore, these masas $A$ possessed the property that $[\Normo{A}{n+1}:\Normo{A}{n}]=2^{\lambda_n}$ for any sequence $(\lambda_n)$ in $\mathbb N$. However, the Jones index is not the best invariant for the inclusion $\Normo{A}{n}\subset\Normo{A}{n+1}$.  In \cite{Jones.PropertiesMasas}, Jones and Popa observe, using work of Ocneanu \cite{Ocneanu.ActionsGroups}, that given an irreducible regular subfactor $S$ of the hyperfinite \IIi factor $R$, the group $G=\Norm{S}/\Unit{S}$ is a countable amenable discrete group, where $\Unit{S}$ is the group of unitaries in $S$.  Furthermore, they show that this $G$ is a complete conjugacy invariant for the inclusion $S\subset R$ and produce, for any countable amenable discrete $G$, a semi-regular masa $A$ in $R$ of length $2$ such that $\Norm{\Normo{A}{1}}/\Unit{\Normo{A}{1}}\cong G$.  

In this paper we present the definitive result in this area, Theorem \ref{Intro.Main} below, that semi-regular masas can be found in the hyperfinite \IIi factor with any countable ordinal length and any collection of countable discrete amenable groups describing the repeated inclusions $\Normo{A}{\beta}\subset\Normo{A}{\beta+1}$.  The case $\alpha=2$ can be found in Theorem 4.1 of \cite{Jones.PropertiesMasas}.
\begin{thm}\label{Intro.Main}
Let $\alpha>1$ be a countable ordinal.  Given countable amenable discrete groups $G_\beta\neq 1$ for $1\leq \beta<\alpha$, there is a semi-regular masa $A$ in the hyperfinite \IIi factor $R$ with length $\alpha$, $\Normo{A}{\alpha}=R$ and
$$
\Norm{\Normo{A}{\beta}}/\Unit{\Normo{A}{\beta}}\cong G_\beta,
$$
for all $1\leq \beta<\alpha$.
\end{thm}
It should be noted that by taking an appropriate free product, one can achieve semi-regular masas in the (interpolated) free-group factors $\VN{\mathbb F_r}$ --- or indeed any factor $R*Q$ where $Q$ is a finite von Neumann algebra --- of length $\alpha$ and with $\Norm{\Normo{A}{\beta}}/\Unit{\Normo{A}{\beta}}\cong G_\beta$ for any countable ordinal $\alpha$ and non-trivial countable amenable discrete groups $(G_\beta)_{\beta<\alpha}$.  In this case $\Normo{A}{\alpha}$ is an irreducible singular subfactor of $\VN{\mathbb F_r}$.  To see this one only needs to note that, if $B$ is a diffuse subalgebra of the \IIi factor $M$ and if $Q$ is any finite von Neumann algebra, then $\NormN{B}{M*Q}=\NormN{B}{M}$, regarded as a subset of $M*Q$.  This fact can be found in \cite[Theorem 2.3]{Sinclair.FreePuk} or \cite{Popa.IPP}.

We establish Theorem \ref{Intro.Main} using a group construction inspired by \cite{Jones.PropertiesMasas} and \cite{Sinclair.StrongSing2}. This distinguishes our approach from that taken in Tauer's original paper \cite{Tauer.Masa} and the earlier work of the first author \cite{Saw.Thesis}. In section \ref{Masa} we construct masas inside a cross product algebra $N\rtimes K$ whose normalisers generate $N\rtimes H$, for some appropriate subgroup $H$ of $K$.  In section \ref{Norm}, we show that the normalisers of $N\rtimes H$ generate $N\rtimes\NormN{H}{K}$, where $\NormN{H}{K}$ is the group normalisers of $H$ in $K$.  We begin in the next section by constructing an inclusion $H\leq K$ corresponding to families $(G_\beta)$ which will yield the required masa.  This construction relies on a repeated iteration of the wreath product. We assemble the proof of Theorem \ref{Intro.Main} at the end of section \ref{Norm}.

\section{Iterating the wreath product}\label{Wreath}
Recall that the wreath product $K\wr G$ of discrete groups $K$ and $G$ is obtained by considering the group $K^G$ of all functions $f:G\rightarrow K$ such that $f(x)=1_K$ for all but finitely many $x\in G$, equipped with pointwise operations inherited from those of $K$.  We then let $G$ act by translation, i.e. $(g\cdot f)(x)=f(g^{-1}x)$ for $f\in K^G$ and $g,x\in G$.  The wreath product $K\wr G$ is the semi-direct product $K^G\rtimes G$ with this action.  It is countable when $K$ and $G$ are countable and amenable when $K$ and $G$ are amenable.  

In this section, we suppose that $\alpha>1$ is a fixed, non-zero countable ordinal and we fix a family $(G_\beta)_{1\leq\beta<\alpha}$ of countable discrete amenable groups, each of which is non-trivial. Given an inclusion $H\leq K$ of groups, we write $\NormN{H}{K}$ for the group of normalisers of $H$ in $K$, that is all those $k\in K$ with $kHk^{-1}=H$. We wish to find a group $H_{\alpha}$ and subgroups $(H_{\beta})_{1\leq\beta<\alpha}$ of $H_\alpha$ with $\NormN{H_{\beta}}{H_\alpha}=H_{\beta +1}$ and $H_{\beta+1}/H_\beta\cong G_\beta$ for each $\beta<\alpha$. Before considering an arbitrary countable ordinal $\alpha$, it is useful to consider the case where $\alpha=3$, i.e., we have only the groups $G_1$ and $G_2$.  In this case we first form $H_3$ in stages.  Define $K_1=\Z$, $K_2=K_1\wr G_1$ and $K_3=K_2\wr G_2$.  Take $H_3=K_3$.  We now work backwards to obtain $H_2$ and $H_1$ with the desired properties.  Let $H_2=K_2^{G_2}$, so that $H_3=H_2\rtimes G_2$ and $H_3/H_2\cong G_2$.  To obtain $H_1$ we decompose $K_2^{G_2}$ into a direct sum $K_2\oplus L_2$, where $L_2$ consists of those finitely supported functions $l:G_2\rightarrow K_2$ which have $l(1_{G_2})=1_{K_2}$.  This gives the identity
$$
H_2=K^{G_2}_2=K_2\oplus L_2=(\Z^{G_1}\rtimes G_1)\oplus L_2.
$$
Now define $H_1=\Z^{G_1}\oplus L_2$, a normal subgroup of $H_2$ with $H_2/H_1\cong G_1$. This gives us the required sequence of groups when $\alpha=3$. We now return to the general case.

Define an increasing family $(K_\beta)_{1\leq\beta\leq\alpha}$ of countable, discrete amenable groups inductively by
$$
K_\beta=\begin{cases}\mathbb Z&\beta=1\\ K_\gamma\wr G_\gamma&\beta=\gamma+1,\quad\gamma>0\\ \cup_{\gamma<\beta}K_\gamma&\beta\textrm{ a limit ordinal}\end{cases},
$$
and write $K=K_\alpha$. To create an increasing family, and hence clarify the last statement in the previous definition, we regard $K_\gamma$ as a subgroup of $K_{\gamma+1}$ by identifying $K_\gamma$ with the functions $f:G_\gamma\rightarrow K_\gamma$ which have $f(x)=1_{K_\gamma}$ for all $x\neq 1_{K_\gamma}$. This is a subgroup of $K_\gamma^{G_\gamma}$ and hence of $K_{\gamma+1}$.  We also define $L_\gamma$ to be the subgroup of $K_{\gamma}^{G_\gamma}$ consisting of all those finitely supported functions $f:G_\gamma\rightarrow K_\gamma$ with $f(1_{K_\gamma})=1_{K_\gamma}$, so that $K_\gamma\oplus L_\gamma=K_\gamma^{G_\gamma}$ as subgroups of $K_{\gamma+1}$.  Henceforth, we regard all the $K_\gamma$ and $L_\gamma$ as subgroups of $K$. 

For each $1\leq\beta<\alpha$, let
\begin{equation}\label{Wreath.DefH}
H_\beta=K_\beta^{G_\beta}\oplus\bigoplus_{\beta<\gamma<\alpha}L_\gamma\leq K,
\end{equation}
and write $H_\alpha=K$. All the subgroups of $K$ involved in the definition of $H_{\beta}$ commute so it makes sense to write the direct sum in (\ref{Wreath.DefH}). Let us note for future use that $H_1$ is an infinite group, as $K_1$ is infinite. Observe that $H_\beta\leq H_{\beta+1}$, as $K_\beta^{G_\beta}\oplus L_{\beta+1}\leq K_{\beta+1}\oplus L_{\beta+1}=K_{\beta+1}^{G_{\beta+1}}$.

We will first show that the $(H_\beta)_{1\leq\beta\leq \alpha}$ form a chain of normalisers inside $K$.  This breaks up into two parts, that $H_{\beta+1}=\NormN{H_\beta}{K}$ for each $\beta<\alpha$, which is Lemma \ref{Wreath.Sucessor} below, and that $H_\beta=\bigcup_{\gamma<\beta}H_\gamma$ when $\beta$ is a limit ordinal (Lemma \ref{Wreath.Limit}).  Finally, in Lemma \ref{Wreath.Cong}, we check that $H_{\beta+1}/H_\beta\cong G_\beta$ for $\beta<\alpha$.

\begin{lem}\label{Wreath.Sucessor}
For every $1\leq\beta<\alpha$, we have $\NormN{H_\beta}{K}=H_{\beta+1}$.
\begin{proof}
If $\beta +1=\alpha$, then $H_\beta=K_\beta^{G_\beta}$ and $H_\alpha=H_\beta\rtimes G_\beta$, so the result holds. We assume $\beta +1<\alpha$. If we can show that 
\begin{equation}\label{Wreath.Sucessor.1}
\NormN{H_\beta}{K}\cap K_{\gamma}=H_{\beta +1}\cap K_{\gamma}
\end{equation}
for all $\gamma>\beta$, then the lemma will follow by taking $\gamma =\alpha$. We will now establish (\ref{Wreath.Sucessor.1}) by induction on $\gamma$.  This splits into three parts.
\begin{enumerate}
\item $\gamma=\beta +1$. In this case $G_\beta$ normalises $K_\beta^{G_\beta}$. As $G_\beta\leq K_{\beta+1}$ and $K_{\beta+1}$ commutes with each $L_\delta$ for $\delta>\beta$, we see that $G_\beta$ normalises $H_\beta$. Hence 
$$
\NormN{H_\beta}{K}\cap K_{\beta+1}=K_{\beta+1}=H_{\beta+1}\cap K_{\beta+1}.
$$

\item $\gamma$ is a limit ordinal. Our inductive hypothesis is  
$$
\NormN{H_\beta}{K}\cap K_{\delta}=H_{\beta+1}\cap K_{\delta}
$$
for $\beta<\delta<\gamma$. By definition, $K_\gamma=\bigcup_{\delta<\gamma}K_\delta=\bigcup_{\beta<\delta<\gamma}K_\delta$ so that 
$$
\NormN{H_\beta}{K}\cap K_\gamma=\bigcup_{\beta<\delta<\gamma} (\NormN{H_\beta}{K} \cap K_\delta)=\bigcup_{\beta<\delta<\gamma}(H_{\beta+1} \cap K_\delta)=H_{\beta+1}\cap K_\gamma.
$$

\item $\gamma>\beta+1$ is a successor ordinal. Let $\gamma=\delta+1$ for some $\delta>\beta$. Then $K_{\gamma}=K_{\delta} ^{G_\delta}\rtimes G_\delta$. Take $(f,g)$ in $K_{\delta} ^{G_\delta}\rtimes G_\delta$, with $g\ne 1_{G_\delta}$. We will demonstrate that $(f,g)$ does not normalise $H_{\beta}$ by producing an element $l$ in $H_{\beta}\cap K_{\gamma}$ which is not conjugated into $H_{\beta}$ by $(f,g)$.

Since $\delta>\beta$, $K_\beta^{G_\beta}$ is a subgroup of  $K_{\delta} ^{G_{\delta}}$. Furthermore, for all $\beta\leq \varrho <\delta$, 
$$
L_{\varrho}\leq L_{\varrho}\oplus K_{\varrho}=K_{\varrho} ^{G_\varrho}\leq K_{\varrho +1}\leq K_{\delta} \leq K_{\delta} ^{G_{\delta}}.
$$
In this way, we see that $K_\gamma$ does not contain $L_{\varrho}$ for $\varrho\geq \gamma$ and contains $L_{\varrho}$ for all other $\varrho$. As a result,
$$
H_{\beta}\cap K_{\gamma}=(K_\beta^{G_\beta}\oplus \bigoplus_{\beta <\varrho<\alpha}L_{\gamma})\cap K_{\gamma}=(K_\beta^{G_\beta}\oplus \bigoplus_{\beta <\varrho<\gamma}L_{\gamma})\cap K_{\gamma}\leq K_\delta^{G_\delta}.
$$
Ergo $H_{\beta}\cap K_{\gamma}=H_{\beta}\cap K_{\delta}^{G_{\delta}}$. Since $K_{\delta} ^{G_\delta}$ is a proper subset of $K_{\gamma}$, we may choose an element $k\in K_{\gamma}$ which is not in $H_{\beta}\cap K_{\gamma}$. Define a function $l:G_{\delta}\rightarrow K_{\delta}$ by setting 
$$
l(x)=\begin{cases}f(1_{G_{\delta}})^{-1}kf(1_{G_{\delta}})&x=g^{-1}\\ 1_{K_{\delta}}&{\rm otherwise} \end{cases}.
$$
By assumption $g\ne 1_{G_{\delta}}$, so that $l(1_{G_{\delta}})=1_{K_{\delta}}$ and $l\in L_{\delta}$. As $\delta>\beta$,  $l\in L_{\delta}\leq H_{\beta}$. If we regard $l$ as $(l,1_{G_{\delta}})$ in $K_{\gamma}=K_\delta^{G_\delta}\rtimes G_\delta$, then
$$
(f,g)(l, 1_{G_{\delta}})(f,g)^{-1}=(f (g\cdot l)f^{-1},1_{G_{\delta}}).
$$

Now for all $x$ in $G_{\delta}$, $(f(g\cdot l)f^{-1})(x)=f(x)l(g^{-1}x)f^{-1}(x)$ which is $1_{K_{\delta}}$ except when $x= 1_{G_{\delta}}$ according to the definition of $l$, and so $f(g\cdot l)f^{-1}$ is in $K_{\delta}$. Recall that in the inclusion of $K_{\delta}$ into $K_{\gamma}$, the element $f(g\cdot l)f^{-1}$ is identified with the element in $(f(g\cdot l)f^{-1})(1_{G_{\delta}})=f(1)l(g^{-1})f(1)^{-1}=k$ in $K_{\delta}$, which was chosen not to lie in $H_{\beta}\cap K_{\gamma}$. Therefore, the element $(f,g)\in K_{\gamma}$ does not normalize $H_{\beta}$, if $g\ne 1_{G_{\delta}}$.

From this, we see that $\NormN{H_{\beta}}{K}\cap K_{\gamma}$ is equal to $\NormN{H_{\beta}}{K}\cap K_{\delta} ^{G_{\delta}}$, which we decompose as $(\NormN{H_{\beta}}{K} \cap K_{\delta})\oplus (\NormN{H_{\beta}}{K} \cap L_{\delta})$. As $L_{\delta}$ is a subgroup of $H_{\beta}$ and hence of $H_{\beta +1}$, $\NormN{H_{\beta}}{K} \cap L_{\delta}=H_{\beta +1} \cap L_{\delta}=L_{\delta}$. Our inductive hypothesis is that $ \NormN{H_{\beta}}{K} \cap K_{\delta}=H_{\beta +1}\cap K_{\delta}$, and so 
$$
\NormN{H_{\beta}}{K}\cap K_{\gamma}=(H_{\beta +1}\cap K_{\delta})\oplus (H_{\beta +1} \cap L_{\delta})=H_{\beta +1}\cap K_{\gamma}^{G_{\gamma}}.
$$

The argument that $H_{\beta}\cap K_{\gamma}=H_{\beta}\cap K_{\delta}^{G_{\delta}}$ also shows $H_{\beta+1}\cap K_{\gamma}=H_{\beta+1}\cap K_{\delta}^{G_{\delta}}$, and this establishes (\ref{Wreath.Sucessor.1}) when $\gamma$ is a sucessor ordinal.
\end{enumerate}
This completes the inductive proof of (\ref{Wreath.Sucessor.1}) and so establishes the lemma.
\end{proof}
\end{lem}

\begin{lem}\label{Wreath.Limit}
Suppose that $\beta\leq\alpha$ is a limit ordinal.  Then 
$$
\bigcup_{\gamma<\beta}H_\gamma=H_\beta.
$$
\begin{proof}
Since the family $(H_\gamma)$ is increasing, we only have to show $\bigcup_{\gamma<\beta}H_\gamma\supset H_\beta$.  Certainly $\bigcup_{\gamma<\beta}H_\gamma\supset K_\beta$, since $K_\beta=\bigcup_{\gamma<\beta}K_\gamma$ and $K_\gamma\leq H_\gamma$.  Since $K_\beta\oplus L_\beta=K_\beta^{G_\beta}$, and $L_\beta\leq H_\gamma$ when $\gamma<\beta$ we have
$$
\bigcup_{\gamma<\beta}H_\gamma\supset K_\beta^{G_\beta}.
$$
For $\beta<\delta<\alpha$, we have $L_\delta\subset\bigcup_{\gamma<\beta}H_\gamma$ and the result follows.
\end{proof}
\end{lem}

\begin{lem}\label{Wreath.Cong}
For each $\beta<\alpha$, $H_{\beta+1}/H_{\beta}\cong G_\beta$.  
\begin{proof}
The group $H_\beta$ is normal in $H_{\beta+1}$ by Lemma \ref{Wreath.Sucessor}. If $\beta+1\neq\alpha$, then
$$
\frac{H_{\beta+1}}{H_\beta}\cong \frac{K_{\beta+1}^{G_{\beta+1}}}{K_\beta^{G_\beta}\oplus L_{\beta+1}}=\frac{K_{\beta+1}\oplus L_{\beta+1}}{K_\beta^{G_\beta}\oplus L_{\beta+1}}\cong\frac{K_{\beta+1}}{K_\beta^{G_\beta}}\cong G_\beta,
$$
where the first isomorphism follows by removing the direct sum $\oplus_{\beta+1<\delta<\alpha}L_\delta$ from both $H_{\beta+1}$ and $H_\beta$, and the last follows from the definition of $K_{\beta+1}$.  If $\beta+1=\alpha$, then $H_\beta=K_\beta^{G_{\beta}}$ and $H_\alpha=K_\alpha=K_\beta\wr G_\beta$, so the result is immediate.
\end{proof}
\end{lem}

\section{A masa in a crossed product factor}\label{Masa}
In this section we produce a masa in a crossed product factor.  These masas will turn out to be the masas whose existance was promised in Theorem \ref{Intro.Main}. The construction is based on Theorem 4.1 of \cite{Jones.PropertiesMasas} but adjusted to ensure the normalisers of the masa generate a subalgebra of the form $N\rtimes H$.

Let $K$ be any countable, infinite, discrete, amenable group and fix an infinite subgroup $H$ of $K$. For each $k\in K$, denote by $R_k$ a copy of the hyperfinite \IIi factor $R$.  Take Cartan masas $A_k\subset R_k$ in each copy of $R$ such that
\begin{enumerate}
\item $A_{k_1}\perp A_{k_2}$, if  $Hk_1\neq Hk_2$;
\item $A_{k_1}=A_{k_2}$, if $Hk_1=Hk_2$,
\end{enumerate}
using the definition of orthogonality of pairs of subalgebras given in \cite{Popa.Orth}.  Such a family can be found; in the proof of Theorem 4.1  of \cite{Jones.PropertiesMasas} an infinite family of pairwise orthogonal masas in the hyperfinite \IIi factor is produced.

Let $N=\vnotimes _{k\in K} R_k$, and consider the Bernoulli shift action $\theta$ of $K$ on $N$ by $\theta_k(\otimes _{g\in K}x_g)=\otimes _{g\in K} x_{kg}$. It is well known that this is a proper outer action, so the crossed product $N\rtimes_{\theta} K$ is a copy of the hyperfinite \IIi factor, which we denote by $M$.  Let $A=\vnotimes_{k\in K}A_k$, a Cartan masa in $N$ which we regard as an abelian subalgebra of $M$.  
\begin{thm}\label{Masa.MainThm}
With the notation above, $A$ is a masa in $M$ and $\NormN{A}{M}''=N\rtimes_{\theta} H$
\begin{proof}
First, we show that $({\mathcal N}_M(A))''=N\rtimes _{\theta} H$. Let $(u_k)_{k\in K}$ denote the implementing unitaries of $\theta$ in $M$. Since $A_{hk}=A_k$, for every $k\in K$ and $h\in H$, the unitaries $(u_h)_{h\in H}$ normalise $A$.  Hence $N\rtimes _{\theta} H \subset \NormN{A}{M}''$. If $k\notin H$, then $A_{kg}\perp A_g$ for all  $g\in K$, and hence, $u_kAu_k ^*\perp A$, by \cite[Lemma 2.4]{Popa.Orth}. Corollary 2.7 of \cite{Popa.Orth} then implies that $u_k\perp\NormN{A}{M}''$.  Since $A$ is Cartan in $N$ we see that $Nu_k$ is orthogonal to $\NormN{A}{M}''$ so that $\NormN{A}{M}''=N\rtimes_\theta H$.

To show that $A$ is a masa, suppose $x\in A'\cap M$.  By the preceeding paragraph, $A'\cap M\subset N\rtimes_\theta H$, so we can write $x=\sum _{h\in H}x_hu_h$ with convergence in strong operator topology, for some $x_h\in N$. Now $ax_h=x_h\theta_h(a)$ for all $a\in A,h\in H$, so $x_{1_H}\in A$. Fix $h\in H$ with $h\neq 1_H$. For $\varepsilon>0$, find a finite set $F$ of $K$ and $y_0\in\vnotimes_{f\in F}R_f$ with $\nm{x_h-y_0}_2<\varepsilon$.  Since $H$ is infinite, there exists $k\in H$ with $k\notin F$ and $k\notin hF$.  Now define $a=\otimes_{g\in K}a_g$ where 
$$
a_g=\begin{cases}1&g\neq k\\u&g=k\end{cases},
$$
and $u$ is a unitary in $A_k$ with $\tr(u)=0$.  Both $a$ and $\theta_h(a)$ are elementary tensors with $1$ in all the positions coming from the finite set $F$. Hence
\begin{equation}\label{Masa.MainThm.1}
\nm{ay_0-y_0\theta_h(a)}_2=\nm{y_0}_2\nm{a-\theta_h(a)}_2=\sqrt{2}\nm{y_0}_2.
\end{equation}
On the other hand, 
\begin{equation}\label{Masa.MainThm.2}
\nm{ay_0-y_0\theta_h(a)}_2\leq 2\nm{a}\nm{x_h-y_0}_2<2\varepsilon.
\end{equation}
Combining equations (\ref{Masa.MainThm.1}) and (\ref{Masa.MainThm.2}) gives
$$
\nm{x_h}_2\leq \nm{x_h-y_0}_2+\nm{y_0}_2<\epsilon+\sqrt{2}\epsilon,
$$
and as $\varepsilon>0$ was arbitrary, $x_h=0$.  Hence $x=x_{1_H}\in A$, so that $A$ is a masa in $M$.
\end{proof}
\end{thm}

\section{Normalisers of cross products}\label{Norm}
In this section our objective is to establish the following theorem regarding the normalisers of cross products by subgroups, which is the last ingredient in the proof of Theorem \ref{Intro.Main}.
\begin{thm}\label{Norm.Cross}
Let $N$ be a \IIi factor and let $\theta$ be a proper outer action of a countable discrete group $K$ on $N$.  For any subgroup $H$ of $K$ we have
$$
\Norm{N\rtimes_\theta H}''=N\rtimes_\theta\NormN{H}{K}.
$$
\end{thm}

Since the intermediate subfactors $N_0$ between $N\rtimes_\theta H$ and $N\rtimes_\theta K$ are all of the form $N\rtimes_\theta G$ for some intermediate group $H\leq G\leq K$ \cite{Choda.Galois}, it is tempting to try and establish Theorem \ref{Norm.Cross} by noting that group elements outside $\NormN{H}{K}$ do not normalise $N\rtimes_\theta H$.  However, this does not exclude the possibility that some $k\in K\setminus\NormN{H}{K}$ can be written as a linear combination of normalisers, so some further argument is needed. 

Throughout this section we let $\theta$ be a proper outer action of the countable discrete group $K$ on the \IIi factor $N$ such as the Bernoulli shift action of the previous section.  Let $N$ be faithfully represented on $\HS$, so that $N\rtimes_\theta K$ is represented on $\HS\otimes_2\ell^2(K)$. We write bounded operators $T$ on $\HS\otimes_2\ell^2(K)$ as matrices $(T_{g,k})_{g,k\in K}$ of elements in $\mathbb B(\HS)$.
\begin{lem}\label{Norm.Lem}
With the notation above, let $\mathcal A\subset\mathbb B(\HS\otimes_2\ell^2(K))$ consist of those operators $T$ with $T_{g,k}=\delta_{g,k}z\theta_k(x)$ for some $x\in N$ and $z\in N'\cap\mathbb B(\HS)$.  Then $\mathcal A''$ consists of the operators $T$ with $T_{g,k}=\delta_{g,k}X_k$ for any uniformly bounded collection of operators $(X_k)_{k\in K}$ in $\mathbb B(\HS)$.
\begin{proof}
Take some $S\in\mathcal A'$ and write $S=(S_{g,k})_{g,k\in K}$.  The form of the operators in $\mathcal A$ ensures that
\begin{equation}\label{Norm.Lem.1}
S_{g,k}z\theta_k(x)=z\theta_g(x)S_{g,k}
\end{equation}
for all $g,k\in K$, $x\in N$ and $z\in N'\cap\mathbb B(\HS)$.  Taking $x=1$ we see that each $S_{g,k}\in N$.  When $g=k$, take $z=1$ so that $S_{g,g}$ lies in $N\cap N'=\mathbb C1$.  Finally, when $g\neq k$, we get $S_{g,k}y=\theta_{gk^{-1}}(y)S_{g,k}$ for all $y\in N$ so that $S_{g,k}=0$, since $\theta$ is a proper outer action.  Ergo $\mathcal A'$ consists of those operators with scalar entries down the diagonal, from which the conclusion readily follows.
\end{proof}
\end{lem}

Theorem \ref{Norm.Cross} follows immediately from the next lemma, which tells us slightly more about the structure of these subfactors.
\begin{lem}\label{Norm.Location}
With the standing notation of this section, let $H$ be a subgroup of $K$.  Suppose that $x\in N\rtimes_\theta K$ satisfies 
$$
x(N\rtimes_\theta H)x^*\subset N\rtimes_\theta H.
$$
Then $x\in N\rtimes_\theta \NormN{H}{K}$.
\begin{proof}
Write $x=\sum_{k\in K}x_ku_k$, where $u_k$ are the unitaries implementing the action by $\theta$ and the sum converges in strong-operator topology. Fix $k\not\in\NormN{H}{K}$ and find $h\in H$ with $khk^{-1}\not\in H$.  For each $g\in K$, let $\phi(g)$ denote the unique element of $K$ with $gh\phi(g)^{-1}=khk^{-1}$.  We will show that $x_k=0$.  

For any $y\in N$
\begin{equation}\label{Norm.Location.1}
\sum_{g\in G}x_g\theta_g(y)\theta_{khk^{-1}}(x_{\phi(g)}^*)=0,
\end{equation}
since this is the $khk^{-1}$ coefficient of $x(yu_h)x^*$ which lies in $N\rtimes_\theta H$ by hypothesis.  Define operators $R$ and $C$ on $\HS\otimes\ell^2(K)$ by
$$
R_{s,t}=\begin{cases}x_t&s=1_K\\0&\textrm{otherwise}\end{cases},\quad C_{s,t}=\begin{cases}\theta_{khk^{-1}}(x_{\phi(s)}^*)&t=1_K\\0&\textrm{otherwise}\end{cases}.
$$

Let $D$ be the diagonal operator with $D_{t,t}=\theta_t(y)$. Then $RDC=0$, since the only entry which does not immediately vanish is the $(1_K,1_K)$ entry, which evaluates to $0$ by (\ref{Norm.Location.1}).  Since $R$ commutes with diagonal matrices with constant diagonal entries $z\in N'$, we see that $R\mathcal AC=0$, where $\mathcal A$ is the algebra described in Lemma \ref{Norm.Lem}.  Hence $R(\mathcal A'')C=0$, and consideration of $T\in\mathcal A''$ with $T_{t,t}=0$ unless $t=k$ gives
$$
x_kT_{k,k}\theta_{khk^{-1}}(x_{k}^*)=0
$$
for all $T_{k,k}\in\mathbb B(\HS)$ from which $x_k=0$ as claimed.
\end{proof}
\end{lem}

We are now able to prove the main result of the paper.
\begin{proof}[Proof of Theorem \ref{Intro.Main}.]
Let $\alpha>1$ be a countable ordinal and $(G_\beta)_{1\leq\beta<\alpha}$ be non-trivial countable discrete amenable groups.  Form the group $K=K_\alpha$ and the infinite subgroup $H=H_1$ described in section \ref{Wreath}. Let $N$ be the infinite tensor product $(\otimes_{k\in K}R_k)''$ of copies of the hyperfinite \IIi factor and $M=N\rtimes_\theta K$, where $\theta$ is the shift action of section \ref{Masa}.  Since $K$ is amenable, $M$ is a copy of the hyperfinite \IIi factor.  Theorem \ref{Masa.MainThm} gives us a masa $A$ in $M$ with 
$$
\NormN{A}{M}''=\Normo{A}{1}=N\rtimes_\theta H_1.
$$
Suppose inductively that we have shown $\Normo{A}{\gamma}=N\rtimes_\theta H_\gamma$ for all $\gamma<\beta$.  If $\beta=\gamma+1$ is a successor ordinal, then Theorem \ref{Norm.Cross} shows that $\Norm{N\rtimes_\theta H_\gamma}''=N\rtimes_\theta\NormN{H_\gamma}{K}$.  Lemma \ref{Wreath.Sucessor} gives $\NormN{H_\gamma}{K}=H_\beta$ so that $\Normo{A}{\beta}=N\rtimes_\theta H_\beta$.  If $\beta$ is a limit ordinal then
$$
\Normo{A}{\beta}=\left(\bigcup_{\gamma<\beta}\Normo{A}{\gamma}\right)''=\left(\bigcup_{\gamma<\beta}N\rtimes_\theta H_\gamma\right)''=N\rtimes_\theta H_\beta,
$$
by Lemma \ref{Wreath.Limit}.  Hence
$$
\Normo{A}{\beta}=N\rtimes_\theta H_\beta,
$$
for all $1\leq\beta\leq\alpha$.  Since $H_\alpha=K$ and $H_\beta\neq K$ for $\beta<\alpha$ the length of $A$ is $\alpha$.  For $\beta<\alpha$,
$$
\Norm{\Normo{A}{\beta+1}}/\Unit{\Normo{A}{\beta}}\cong H_{\beta+1}/H_\beta\cong G_\beta,
$$
where the last identity is Lemma \ref{Wreath.Cong}.
\end{proof}  

\section{The Tauer property}\label{Infinite}
Recall that in section \ref{Intro} we said that an inclusion $N\subset M$ of \IIi factors has the Tauer property if every $x,y\in M$ with $\ce{N}{x}=\ce{N}{y}=0$ satisfy $xy^*\in N$.  Here the map $\mathbb E_N$ is the unique trace-preserving conditional expectation from $M$ onto $N$.  This property was used by Tauer to distinguish between subfactors and hence between masas in the hyperfinite \IIi factor, \cite[Section 6]{Tauer.Masa}.  In this section we establish the next result showing that the Tauer property can be characterised in terms of the Jones index.
\begin{thm}\label{Intro.TauerProp}
An inclusion of subfactors $N\subset M$ has the Tauer property if and only if $[M:N]=2$.
\end{thm}

With the additional assumption that $[M:N]<\infty$ this characterisation is straight-forward. 
\begin{proof}[Proof of Theorem \ref{Intro.TauerProp} for finite index inclusions.]
Using Goldman's Theorem \cite{Goldman.Index2} (see also \cite[Corollary 3.4.3]{Jones.Index}), it is easy to see that any index $2$ inclusion of \IIi factors has the Tauer property.  For the reverse direction, first suppose that $[N:M]<\infty$ and that $N\subset M$ has the Tauer property.  By \cite[Lemma 3.1.8]{Jones.Index}, there is some subfactor $P$ of $N$ with $[N:P]=[M:N]$ such that $M$ is isomorphic to the basic construction $\ip{N}{e_P}$.  Under this identification recall that $\ce{N}{e_P}=[M:N]^{-1}1$.  Take $x=y=1-[M:N]e_P$ so that $\ce{N}{x}=\ce{N}{y}=0$.  Then $xy=1-2[M:N]e_P+[M:N]^2e_P$, lies in $N$ by the Tauer property, so $2[M:N]=[M:N]^2$ and $[M:N]=2$.  
\end{proof}

Establishing that no infinite index inclusion of subfactors can have the Tauer property is more technical.  We use a Pimsner-Popa orthogonal basis argument; the technicalities arise as we are forced to use an unbounded basis to describe a general infinite index inclusion.

Let $N\subset M$ be an inclusion of \IIi factors.  As is usual we assume $M$ is in standard form, that is represented as operators of left-multiplication on $L^2(M)$, the completion of $M$ in $\nm{x}_2=\tr(x^*x)^{1/2}$, where $\tr$ is the unique trace on $M$.  We also complete $M$ in the $1$-norm given by 
$$
\nm{x}_1=\sup_{\substack{y\in M,\nm{y}\leq 1}}|\tr(xy)|
$$ 
and regard $L^2(N)$ and $L^1(N)$ as the closures of $N$ in $L^2(M)$ and $L^1(M)$ respectively.  We have $\nm{x^*y}_1\leq\nm{x}_2\nm{y}_2$ for $x,y\in M$ by Cauchy-Schwartz.  We can use this to define an element $\xi^*\eta\in L^1(M)$ for any $\xi,\eta\in L^2(M)$ by taking sequences $(x_n)$ and $(y_m)$ in $M$ converging to $\xi$ and $\eta$ in $L^2(M)$ respectively.  Then $\xi^*\eta$ is the $L^1$-limit of $x_n^*y_m$ as $n,m\rightarrow\infty$ and is independent of the choice of these sequences. Furthermore
\begin{equation}\label{Infinite.1Norm}
\nm{\xi^*\eta}_1=\sup_{\substack{x,y\in M\\\nm{x},\nm{y}\leq 1}}\ip{Jx^*J\xi}{Jy^*J\eta}\leq\nm{\xi}_2\nm{\eta}_2,
\end{equation}
where $J$ is the conjugation operator on $L^2(M)$ obtained by extending $x\mapsto x^*$.  The orthogonal projection $e_N:L^2(M)\rightarrow L^2(N)$ restricts to the conditional expectation $\mathbb E_N:M\rightarrow N$, which in turn extends by continuity to a contraction $E_N$ of $L^1(M)$ onto $L^1(N)$.  The following reformulation of the Tauer property follows by approximating elements in $L^2(M)$ with elements in $M$.
\begin{prop}\label{Infinite.Tauer}
Let $N\subset M$ be an inclusion of \IIi factors with the Tauer property.  Given any $\xi,\eta\in L^2(M)$ with $e_N(\xi)=e_N(\eta)=0$ we have $\xi^*\eta\in L^1(N)$, i.e. $E_N(\xi^*\eta)=\xi^*\eta$.
\end{prop}

Now suppose that $[M:N]=\infty$.  The basic construction $\ip{M}{e_N}$ of \cite{Jones.Index} is a $\mathrm{II}_\infty$ factor, so we can find unitaries $(u_n)_{n=1}^\infty$ in $\ip{M}{e_N}$ with
$$
\sum_{n=1}^\infty u_ne_Nu_n^*=1-e_N.
$$
Define $\xi_n=u_n\hat{1}\in L^2(M)$, where $\hat{1}$ is the identity of $M$ thought of as a vector in $L^2(M)$.  Since $u_ne_Nu_n^*$ is a projection underneath $1-e_N$ we have $e_Nu_ne_N=0$.  For $x,y\in N$ we have
\begin{align*}
&\ip{Jx^*J\xi_n}{Jy^*J\hat{1}}=\ip{Jx^*Ju_n\hat{1}}{y\hat{1}}=\ip{Jx^*Ju_ne_N\hat{1}}{e_Ny\hat{1}}\\
=&\ip{e_NJx^*Ju_ne_N\hat{1}}{y\hat{1}}=\ip{Jx^*Je_Nu_ne_N\hat{1}}{y\hat{1}}=0,
\end{align*}
as $Jx^*J$ and $Jy^*J$ commute with $e_N$.  Hence $e_N(\xi_n)=0$.  Similarly, for $m\neq n$ we have $e_Nu_m^*u_ne_N=0$ and this gives $E_N(\xi_n^*\xi_m)=0$ as
\begin{align*}
&\ip{Jx^*J\xi_n}{Jy^*J\xi_m}=\ip{Jx^*Ju_n\hat{1}}{Jy^*Ju_m\hat{1}}=\ip{Jx^*Ju_ne_N\hat{1}}{Jy^*Ju_me_N\hat{1}}\\
=&\ip{u_ne_NJx^*J\hat{1}}{u_me_NJy^*J\hat{1}}=\ip{Jx^*J\hat{1}}{e_Nu_n^*u_me_NJy^*J\hat{1}}=0,
\end{align*}
for all $x,y\in N$.

By way of obtaining a contradiction, suppose that $N\subset M$ has the Tauer property.  Proposition \ref{Infinite.Tauer} gives $\xi_n^*\xi_m\in L^1(N)$ for all $m,n$ since each $e_N(\xi_n)=0$.  For $m\neq n$, we have $\xi_n^*\xi_m=E_N(\xi_n^*\xi_m)=0$.  On the other hand $u_m$ and $u_n$ are unitaries in $\ip{M}{e_N}$ so $u_nu_m^*\neq 0$.  By density of $M$ in $L^2(M)$ there are $x,y\in M$ with $\ip{u_nx}{u_my}\neq0$.  Then
$$
\ip{u_nx}{un_y}=\ip{Jx^*Ju_n\hat{1}}{Jy^*Ju_m\hat{1}}=\ip{Jx^*J\xi_n}{Jy^*J\xi_m}\neq 0
$$
so that $\xi_n^*\xi_m\neq 0$ by (\ref{Infinite.1Norm}).  This gives the required contradiction, so that no infinite index inclusion of \IIi factors can have the Tauer property, and thus completes the proof of Theorem \ref{Intro.TauerProp}.

\begin{ack}
The authors would like to thank Roger Smith for his helpful suggestions and conversations regarding the work in this paper.  The first author would like to thank Allan Sinclair for suggesting the notion of ordinal lengths.
\end{ack}

\bigskip

Department of Mathematics, Texas A\&M University, College Station, TX77843, U.S.A.
\texttt{white@math.tamu.edu, awiggins@math.tamu.edu}
\end{document}